\documentclass{agtart_a}
\pdfoutput=1

\usepackage{amscd,enumerate,verbatim}
\usepackage[all]{xy}

\title{Aspherical manifolds, relative hyperbolicity, simplicial volume and
assembly maps}

\author{Igor Belegradek}
\givenname{Igor}
\surname{Belegradek}
\address{School of Mathematics\\
Georgia Institute of Technology\\\newline
Atlanta GA 30332--0160}
\email{ib@math.gatech.edu}
\urladdr{}

\volumenumber{6}
\issuenumber{}
\publicationyear{2006}
\papernumber{48}
\startpage{1341}
\endpage{1354}

\doi{}
\MR{}
\Zbl{}


\received{19 October 2005}
\revised{}
\accepted{30 June 2006}
\published{20 September 2006}
\publishedonline{20 September 2006}
\proposed{}
\seconded{}
\corresponding{}
\editor{}
\version{}

\arxivreference{math.GR/0509504}

\let\xysavmatrix\xymatrix
\def\xymatrix{\disablesubscriptcorrection\xysavmatrix}
\AtBeginDocument{\let\bar\wbar\let\tilde\wtilde\let\hat\what}

\makeop{curv}
\makeop{inj}
\makeop{vol}
\makeop{diam}
\makeop{Ric}
\makeop{Iso}
\makeop{Fix}
\makeop{Int}
\makeop{pinch}
\makeop{cd}

\def\d{\partial}

\def\R{\mathbb{R}}
\def\Z{\mathbb{Z}}

\makeatletter
\def\cnewtheorem#1[#2]#3{\newtheorem{#1}{#3}[section]
\expandafter\let\csname c@#1\endcsname\c@thm}
\makeatother

\newtheorem{thm}{Theorem}[section]
\cnewtheorem{cor}[thm]{Corollary}
\cnewtheorem{quest}[thm]{Question}
\cnewtheorem{lem}[thm]{Lemma}
\cnewtheorem{prop}[thm]{Proposition}
\cnewtheorem{defn}[thm]{Definition}
\cnewtheorem{conj}[thm]{Conjecture}
\cnewtheorem{Example}[thm]{Example}
\makeautorefname{Example}{Example}
\newenvironment{ex}{\begin{Example}\rm}{\end{Example}}
\cnewtheorem{remark}[thm]{Remark}
\newenvironment{rmk}{\begin{remark}\rm}{\end{remark}}
\cnewtheorem{Fact}[thm]{Fact}

\cnewtheorem{Main Lemma}[thm]{Main Lemma}

\cnewtheorem{Convention}[thm]{Convention}

\newtheorem{amc}{Assembly Map Conjecture}

\makeautorefname{amc}{Assembly Map Conjecture}
\newtheorem{nhsc}{Novikov Higher Signatures Conjecture}

\makeautorefname{nhsc}{Novikov Higher Signatures Conjecture}
\newtheorem{rbc}{Relative Borel Conjecture}

\makeautorefname{rbc}{Relative Borel Conjecture}

\begin{document}

\begin{asciiabstract}
This paper contains examples of closed aspherical manifolds obtained
as a by-product of recent work by the author [arXiv:math.GR/0509490] on
the relative strict hyperbolization of polyhedra. The following is proved.

(I) Any closed aspherical triangulated n-manifold M^n with
hyperbolic fundamental group is a retract of a closed aspherical
triangulated (n+1)-manifold N^(n+1) with hyperbolic fundamental
group.

(II) If B_1,...,B_m are closed aspherical triangulated
n-manifolds, then there is a closed aspherical triangulated manifold
N of dimension n+1 such that N has nonzero simplicial volume,
N retracts to each B_k, and \pi_1(N) is hyperbolic relative to
\pi_1(B_k)'s.

(III) Any finite aspherical simplicial  complex is a retract of a
closed aspherical triangulated manifold with positive simplicial volume
and non-elementary relatively hyperbolic fundamental group.
\end{asciiabstract}

\begin{htmlabstract}
<p class="noindent">
This paper contains examples of closed aspherical manifolds obtained as
a by-product of recent work by the author [arXiv:math.GR/0509490] on the relative
strict hyperbolization of polyhedra. The following is proved.
</p>
<p class="noindent">
(I) Any closed aspherical triangulated n&ndash;manifold M<sup>n</sup> with
hyperbolic fundamental group is a retract of a closed aspherical
triangulated (n+1)&ndash;manifold N<sup>n+1</sup> with hyperbolic fundamental
group.
</p>
<p class="noindent">
(II) If B<sub>1</sub>,&hellip; B<sub>m</sub> are closed aspherical triangulated
n&ndash;manifolds, then there is a closed aspherical triangulated manifold
N of dimension n+1 such that N has nonzero simplicial volume,
N retracts to each B<sub>k</sub>, and &pi;<sub>1</sub>(N) is hyperbolic relative to
&pi;<sub>1</sub>(B<sub>k</sub>)'s.
</p>
<p class="noindent">
(III) Any finite aspherical simplicial  complex is a retract of a
closed aspherical triangulated manifold with positive simplicial volume
and non-elementary relatively hyperbolic fundamental group.
</p>
\end{htmlabstract}

\begin{abstract}
This paper contains examples of closed aspherical manifolds obtained as
a by-product of recent work by the author~\cite{Bel-hp} on the relative
strict hyperbolization of polyhedra. The following is proved.

(I)\hskip .75em Any closed aspherical triangulated $n$--manifold $M^n$ with
hyperbolic fundamental group is a retract of a closed aspherical
triangulated $(n{+}1)$--manifold $N^{n+1}$ with hyperbolic fundamental
group.

(II)\hskip .75em If $B_1,\dots B_m$ are closed aspherical triangulated
$n$--manifolds, then there is a closed aspherical triangulated manifold
$N$ of dimension $n{+}1$ such that $N$ has nonzero simplicial volume,
$N$ retracts to each $B_k$, and $\pi_1(N)$ is hyperbolic relative to
$\pi_1(B_k)$'s.

(III)\hskip .75em Any finite aspherical simplicial  complex is a retract of a
closed aspherical triangulated manifold with positive simplicial volume
and non-elementary relatively hyperbolic fundamental group.
\end{abstract}

\maketitle

\section{Introduction}
\label{sec: intro}

Recently there has been a surge of interest in relatively 
hyperbolic groups (see Drutu and Sapir~\cite{DruSap}, Osin~\cite{Osi-rel},
Dahmani~\cite{Dah-acc}, Groves~\cite{Groves} and Yaman~\cite{Yaman}). 
Three basic classes of examples are finite free products 
(which are hyperbolic relative to the factors), hyperbolic groups
(which are hyperbolic relative to the trivial subgroups),
and geometrically finite isometry groups of negatively
pinched Hadamard manifolds (which are hyperbolic relative
to the maximal parabolic subgroups). Finally,
any group is hyperbolic relative to itself; we call 
a relatively hyperbolic group {\it non-elementary} 
unless it is finite, or virtually--$\Z$, or hyperbolic relative 
to itself.

A major source of examples of relatively hyperbolic groups is 
the small cancellation theory (see Osin~\cite{Osi-sc}) which typically involves 
a $2$--dimensional construction, such as 
adding a ``sufficiently long'' relation to the group.
In~\cite{Bel-hp} the author showed that the spaces obtained via 
relative strict hyperbolization of polyhedra have relatively 
hyperbolic fundamental groups, which allows to
construct many higher-dimensional relatively 
hyperbolic groups. Here is the main result of this note.

\begin{thm}\label{thm: intro-main}
If $B^n$ is a closed aspherical triangulated $n$--manifold such
that $\pi_1(B^n)$ is hyperbolic relative to the subgroups
$H_1,\dots ,H_k$, then there is a closed aspherical triangulated  
$(n{+}1)$--manifold $N^{n+1}$ and an embedding
$i\co B^n\to N^{n+1}$ such that
\begin{itemize}
\item $N^{n+1}$ retracts onto $i(B^n)$,
\item $N^{n+1}$ has nonzero simplicial volume, and
\item $\pi_1(N^{n+1})$ is hyperbolic relative to 
$i_\ast(H_1),\dots ,i_\ast (H_k)$.
\end{itemize}
\end{thm}

\fullref{thm: intro-main} is also true
in PL and smooth categories, meaning that 
the embedding $i\co B^n\to N^{n+1}$ 
and the retraction $N^{n+1}\to i(B^n)$ are
either PL or smooth. I do not know whether 
\fullref{thm: intro-main} holds in the 
topological category. 
If each $H_i$ is trivial, then 
\fullref{thm: intro-main} reduces to (I), while
(II) follows from \fullref{thm: main}, 
which is a slight generalization of \fullref{thm: intro-main}. 
The item (III) follows from (II) combined with the reflection group 
trick of Davis~\cite{Dav-ann} implying  that 
each finite aspherical simplicial complex $K$
is a retract of a closed aspherical manifold $M$. 

The reflection group trick  was used with great 
success by Davis and collaborators
to produce closed aspherical manifolds
with various exotic properties (see Davis~\cite{Davis-ex} and 
references therein), eg if $\pi_1(K)$ is not residually finite,
then neither is $\pi_1(M)$, because residual finiteness is inherited
by subgroups. 
Similarly, \fullref{thm: intro-main} gives a new source  
of closed aspherical manifolds with nonzero simplicial volume
and various exotic properties.
Adapting an example of Weinberger, we build in each dimension $\ge 5$
a closed aspherical manifold $N$ with nonzero simplicial volume 
such that $\pi_1(N)$ is a non-elementary relatively hyperbolic group 
with unsolvable word problem. 
In particular, $\pi_1(N)$ is not hyperbolic, not $CAT(0)$, 
not automatic, not residually finite, and not linear over 
any commutative ring.

Also (I) yields numerous closed aspherical triangulated manifolds with 
hyperbolic fundamental group, which admit no obvious locally 
$CAT(-1)$--metric. In fact, the hyperbolicity of the fundamental group
follows by an indirect argument that uses Dahmani's combination theorem
for relatively hyperbolic groups~\cite{Dah-comb}, and Osin's result that if
all the maximal parabolic subgroups of a relatively hyperbolic
group have linear Dehn function, then the ambient relatively hyperbolic
group is hyperbolic~\cite{Osi-rel}. 

More generally, if we specify a class  of groups $\mathcal C$, 
eg abelian, nilpotent, slender etc, then given a 
closed aspherical triangulated  manifold $B$ with $\pi_1(B)$ 
hyperbolic relative to subgroups that belong to $\mathcal C$,
\fullref{thm: intro-main} produces manifolds with 
the same properties in each dimension $>\dim(B)$.

Just like in~\cite[Section~11]{Davis-ex} we have
applications to the assembly maps, eg 
(I) implies that if the \fullref{amc}
(that the assembly maps in $L$--theory and $K$--theory 
are isomorphisms) holds for all closed aspherical 
triangulated $(n{+}1)$--manifolds with hyperbolic fundamental groups, 
then it also holds for all closed aspherical triangulated
$n$--manifolds with hyperbolic fundamental groups. 
Similarly, (II) shows that the \fullref{amc}
for closed aspherical triangulated manifolds is equivalent
to the \fullref{amc}
for closed aspherical triangulated 
manifolds with nonzero simplicial volume
and non-elementary relatively hyperbolic fundamental group.
Finally,
(III) implies that if the Novikov Conjecture holds for all 
non-elementary relatively hyperbolic fundamental groups
which are fundamental groups of closed triangulated
aspherical manifold with nonzero simplicial volume, then
it holds for the fundamental groups of all
finite aspherical complexes. All these results on assembly maps
also hold for PL or smooth manifolds.

\section{Relative strict hyperbolization}

We refer to Gromov~\cite{Gro-hypgr}, and to the work of Charney, Davis,
Januszkiewicz and Weinberger ~\cite{CD,DJW,DJ} for general background
on hyperbolization of polyhedra, and to~\cite{Bel-hp} for a detailed
study of relative strict hyperbolization.

For the purposes of this paper the relative strict hyperbolization 
is a procedure that takes as an input a compact  
connected triangulated 
$n$--manifold pair $(M,\d M)$ and produces another 
triangulated $n$--manifold pair $(R,\d R)$ together with a surjective 
continuous map $h\co (R,\d R)\to (M,\d M)$ 
that restricts to a homeomorphism
of the boundaries $\d R\to\d M$, 
induces surjections on homology and the fundamental
group, and pulls back rational Pontrjagin classes and the first 
Stiefel--Whitney class. 
All this also works in smooth category by working with
smooth triangulations. 

Also $\d R$ is incompressible in $R$ (ie no homotopically
nontrivial loop in $\d R$ is null-homotopic in $R$), and
if each path--component of $\d M$ is aspherical, 
then $R$ is aspherical.
Furthermore, the space obtained from $R$ by attaching a cone on $\d R$
is a locally $CAT(-1)$ piecewise-hyperbolic simplicial complex, and
according to~\cite{Bel-hp} $\pi_1(R)$ is non-elementary
relatively hyperbolic, relative to the fundamental groups of
the path-components of $\d R$ that are infinite.

\section{Main result}

All the proofs in the section hold without change
in either for triangulated, PL, or smooth manifolds;
for brevity we just deal with the PL case.

\begin{thm}\label{thm: main} 
If $B_1^n,\dots B^n_m$ are closed aspherical 
PL $n$--manifolds with $n>0$ such
that for each $k=1,\dots m$, the group $\pi_1(B_k^n)$ 
is hyperbolic relative to the subgroups 
$H^k_1,\dots ,H^k_{m_k}$, then there is a closed aspherical PL 
$(n{+}1)$--manifold $N^{n+1}$ and PL embeddings
$i_k\co B_k^n\to N^{n+1}$ such that the following hold.
\begin{itemize}
\item For each $k=1,\dots m$, there is a retraction 
$r_k\co N^{n+1}\to i_k(B_k^n)$.
\item $\pi_1(N^{n+1})$ is hyperbolic relative to 
$\{i_{k\ast}(H^k_j): k=1,\dots m,\ j=1,\dots m_k\}$.
\item $N^{n+1}$ has nonzero simplicial volume.
\end{itemize}
\end{thm}
\begin{proof}
Let $M$ be the connected sum of all the manifolds $B_k\times [0,1]$;
since $n>0$, the manifold $M$ is connected.
For any fixed $k$, we can collapse to the point each summand 
$B_j\times [0,1]$ with $j\neq k$, which defines a retraction 
$M\to B_k\times [0,1]$. Therefore, $M$ retracts onto each $B_k$.
Now apply the relative strict hyperbolization to $(M,\d M)$
to obtain $(R,\d R)$, and look at the double $DR$ of $R$ along $\d R$.
The double $N^{n+1}:=DR$ retracts onto each $B_k$ as follows.
First apply the quotient map of the double
involution, then use $h\co (R,\d R)\to (M,\d M)$, 
then use the retraction $M\to B_k$ constructed above, 
and finally go back to $\d R$ via $h^{-1}$.

By~\cite{Bel-hp}
$\pi_1(R)$ is hyperbolic relative to $\pi_1(B_k)$'s,
so a combination theorem for relatively hyperbolic groups due to 
Dahmani~\cite{Dah-comb} implies that 
$\pi_1(DR)$ is hyperbolic relative to $\pi_1(B_k)$'s.
Hence, by a recent result of Drutu--Sapir~\cite[Corollary 1.14]{DruSap}
$\pi_1(DR)$ is hyperbolic relative to the images of
$H_j^k$ under $\pi_1$--injective inclusions $B_k\to DR$. 

To see that
the double $DR$ has nonzero simplicial volume, note that
the simplicial volume is nonincreasing under continuous maps, and
the quotient map $DR\to DR/R=R/\d R$ maps the fundamental class $[DR]$
onto the the fundamental class of the pseudomanifold $R/\d R$,
which is locally $CAT(-1)$~\cite[Remark 3.2]{Bel-hp};
in fact, $R/\d R$ is the strict hyperbolization
of the space obtained from $M$ by attaching a cone over $\d M$.
Hence $R/\d R$ has positive simplicial volume (see Yamaguchi~\cite{Yamag}).
(Instead of using~\cite{Yamag}, one could employ a more general 
result of Mineyev~\cite{Min} that since $R/\d R$ is aspherical
and $\pi_1(R/\d R)$ is hyperbolic, 
every cohomology class in $H^{n+1}(R/\d R,\R)$ is bounded for $n>0$, 
which implies by a standard argument~ (see Benedetti and
Petronio~\cite[Proposition F.2.2]{Ben-Pet}) that $||R/\d R||>0$).
\end{proof}

\begin{proof}[Proof of (I) of the abstract]
We keep the notations of \fullref{thm: main} with $m=1$.
By assumption $\pi_1(B_1)$ is hyperbolic, and hence, 
by a theorem of Osin~\cite[Corollary 2.41]{Osi-rel} 
%
%
we can remove
$\pi_1(B_1)$'s from the list of maximal parabolic subgroups of
$\pi_1(DR)$, so that $\pi_1(DR)$ is hyperbolic. Alternatively,
instead of referring to the Osin's result one could a result 
of Drutu--Sapir~\cite[Corollary 1.14]{DruSap} to conclude that $\pi_1(DR)$
is hyperbolic relative to the trivial subgroup, and 
hence is hyperbolic.
\end{proof}

\begin{proof}[Proof of (II) of the abstract]
Since any group is hyperbolic relative to itself, we can
let $m_k=1$ and $H^k_1:=\pi_1(B_k)$ for each $k=1,\dots, m$,
so \fullref{thm: main} yields the desired assertion.
\end{proof}

\begin{rmk}
One can show that $DR$ is orientable if and only if each $B_k$ is orientable.
Indeed, if $DR$ is orientable, then so is each $B_k$ because it sits
in $DR$ with trivial normal bundle. Conversely, if each $B_k$ is orientable,
then so is $M$, and hence $R$ is orientable because 
the relative strict hyperbolization preserves orientability. 
Now $DR$ is orientable as the double of an orientable manifold.
\end{rmk}

\begin{rmk}\label{prop: retract w1}
If $m=1$, one can improve the previous remark to the claim that
the retraction $r_1$ pulls back the first Stiefel--Whitney class $w_1$.
In fact, $r_1$ is the composition of three retractions, and 
each of them pulls back $w_1$. This is obvious for 
the projection $M=B_1\times [0,1]\to B_1$ that in fact
pulls back the stable tangent bundle (and this is what fails if $m>1$).
Also the hyperbolization map $R\to M$ pulls back 
$w_1$ essentially because the hyperbolic manifold with corners
used as the building block in the strict hyperbolization is orientable,
see~\cite[Remark 3.2]{Bel-hp}. 
Finally, an elementary covering space argument shows that
the retraction $DR\to R$ always pulls back $w_1$.
\end{rmk}

\section{Closed aspherical manifolds with nonzero simplicial volume}

There are very few known ways to produce closed aspherical manifolds
of nonzero simplicial volume. Denote by $\mathcal {P}$
the class of closed manifolds
of nonzero  simplicial volume, and by 
$\mathcal {AP}\subset\mathcal {P}$ the subclass
of manifolds in $\mathcal P$ that are aspherical. 

The classes $\mathcal P$ and $\mathcal {AP}$
are closed under the following operations:
\begin{enumerate}
\item
products (see Gromov~\cite[page 10]{Grom-vol-bound-coh}).
\item
finite covers, quotients by free actions of
finite groups, and homotopy equivalences to a closed 
manifold~\cite[page 8]{Grom-vol-bound-coh}.
\item
taking surface bundles 
(ie if $V$ is in $\mathcal P$, and
$\Sigma$ is an orientable closed surface of $\chi(\Sigma)<0$, 
then the total space of any smooth orientable $\Sigma$--bundle
over $V$ is in $\mathcal P$, and by the homotopy exact sequence
of the fibration, if the same holds when $\mathcal P$ is replaced
with $\mathcal {AP}$) (see Hoster and Kotschick~\cite{HK}).
\end{enumerate}

\begin{rmk}
\label{rmk: gluing} 
Manifolds in $\mathcal {P}$ can be also produced by
gluing compact manifolds with 
positive relative simplicial volume along incompressible 
boundary components that have amenable fundamental 
groups (see Gromov~\cite[page 55]{Grom-vol-bound-coh}, and
Kuessner~\cite{Kue-mul}.  Furthermore, the resulting manifold is in
$\mathcal {AP}$ if all the pieces are aspherical, and all the boundary
components are aspherical and incompressible (because any graph of
spaces with aspherical edges and vertices is aspherical provided the
edge-to-vertex maps are $\pi_1$--injective).
\end{rmk}

The class $\mathcal {AP}$ contains the following manifolds:

\begin{itemize}
\item
closed negatively curved Riemannian manifolds 
(this is due to Thurston, who also proved positivity of relative
simplicial volume for open complete negatively pinched manifolds 
of finite volume, see~\cite[Section 1.2]{Grom-vol-bound-coh}, 
\item
closed aspherical manifolds 
with hyperbolic fundamental groups (see Mineyev \cite{Min})
for example, the manifolds obtained by strict hyperbolization (see
Charney and Davis~\cite{CD}),
\item
closed locally symmetric manifolds of nonpositive 
curvature (see Lafont and Schmidt~\cite{Laf-Sch}),
\item
doubles of finite volume negatively pinched manifolds. 
More generally,
one can glue several manifolds as in \fullref{rmk: gluing},
provided at least one piece has nonzero relative simplicial volume,
as happens eg
for closed aspherical Haken $3$--manifolds with at least
one hyperbolic piece in the JSJ--decomposition (see Soma~\cite{Soma}).
\end{itemize}

The above lists all sources of manifolds $\mathcal{AP}$ known 
to the author before writing this note. 
There are some other manifolds in $\mathcal{P}$,
eg the manifolds that appear as the base
of flat affine bundle with nonzero Euler 
number (see Gromov~\cite[page 23]{Grom-vol-bound-coh}, Smillie~\cite{Smi}
and Hausmann~\cite{Hau}). 
Unfortunately, I do not know whether any of these examples 
lie in $\mathcal{AP}$, with the obvious exception
of closed orientable surfaces of negative Euler 
characteristic, or their products, which are already
included in the above lists.

In this section we add 
more items to the above list by exploring \fullref{rmk: gluing}.
We start from the following observation.

\begin{prop}\label{thm: closed siml vol} 
If $M$ is a compact orientable PL $n$--manifold, with $n\ge 2$, 
such that the pseudomanifold $M/\d M$ admits a locally--$CAT(-1)$ 
metric. Then any closed orientable
manifold $N$ containing $M$ as a codimension zero PL
submanifold has nonzero simplicial volume.
\end{prop}
\begin{proof}
If $N$ is an arbitrary closed orientable
manifold $N$ containing $M$ as a codimension zero PL--submanifold,
then the quotient $\bar Z$ of $N$, obtained by collapsing 
$N\setminus\mathrm{Int}(M)$ to the point, 
is homeomorphic to $M/\d M$, which is 
an oriented pseudomanifold whose fundamental class is the image
of the fundamental class of $N$ under 
the quotient map $N\to \bar Z$.
Since simplicial volume is nonincreasing under continuous maps, 
we deduce $||N||\ge ||\bar Z||$. The 
cohomology class dual to the fundamental class of $Z$,
is bounded by Mineyev~\cite{Min} as $n\ge 2$. Hence $||\bar Z||>0$
by a standard argument (see, for example, Benedetti and
Petronio~\cite[Proposition F.2.2]{Ben-Pet}).
\end{proof}

\begin{ex}
If $M$ is any manifold of dimension $\ge 2$
obtained obtained by relative strict hyperbolization,
then $M/\d M$ admits a piecewise-hyperbolic locally--$CAT(-1)$ 
metric as was noted in~\cite{Bel-hp}. Note that $M$ is aspherical
if and only if each component of $\d M$ is aspherical; 
in this case the double of $M$ is also aspherical so 
by \fullref{thm: closed siml vol} 
it lies in
$\mathcal{AP}$, and more generally we could glue $M$
with several other aspherical manifolds and the result
will be in $\mathcal{AP}$ provided all their boundary
components are aspherical and incompressible.
\end{ex}

\begin{ex}\label{ex: mosher-sageev}
If $M$ is a compact manifold with boundary
obtained by chopping off cusps of a complete finite volume real 
hyperbolic manifold, and if each component of $\d M$
is a flat torus of injectivity radius $>\pi$,
then  $M/\d M$ admits a locally--$CAT(-1)$ metric (see Mosher and
Sageev~\cite{MS}). 
Again the double of $M$ lies in $\mathcal{AP}$.
Furthermore, there are other ways to embed $M$ into a closed
aspherical manifold, eg by cusp closing (see Schroeder~\cite{Schr}), 
or by Dehn surgery on $\d M$ (see~\cite{And-dehn}, where
Anderson constructs Einstein metrics on these manifolds).
Here the ambient closed manifolds are aspherical, 
because they can be shown to
admit metrics of nonpositive sectional curvature
via Gromov--Thurston $2\pi$--Theorem. 
\end{ex}

\begin{rmk} The assumption of \fullref{ex: mosher-sageev}
that the injectivity radii of flat components of $\d M$ are $>\pi$
cannot be dropped, eg any hyperbolic
knot complement in the $3$--sphere, or the punctured torus
sitting inside the $2$--torus are codimension zero submanifolds
of the closed manifolds whose simplicial volume vanishes. 
\end{rmk}

\begin{rmk} 
Since the fundamental groups of finite volume
complete hyperbolic manifolds are residually finite,
any such manifold has a finite cover satisfying the assumptions
of \fullref{ex: mosher-sageev}.
\end{rmk}

\begin{rmk}\label{gromovs question}
Gromov stated in~\cite[page 189]{Gro-asy-inv} 
that a result similar to \fullref{ex: mosher-sageev}
should be true in all negatively curved symmetric spaces.
I do not know how to prove this claim.
\end{rmk}

\section{Simplicial volume vs nonpositive curvature}

One may naively suspect that all manifolds in $\mathcal{AP}$ 
admit metrics of nonpositive sectional curvature. It seems
this problem was first discussed by Hoster and Kotschick~\cite{HK}, 
who noted that the surface bundles over surfaces 
belong to $\mathcal{AP}$ while some of them admit no
nonpositively curved metrics, by the work of Kapovich and Leeb~\cite{KapLe}.
Also there are manifolds in $\mathcal{AP}$
that are homeomorphic but not diffeomorphic to
nonpositively curved ones (see Aravinda and Farrell~\cite{AF},
Okun~\cite{Okun}, and Ontaneda~\cite{Ont}).

Perhaps, the simplest example can be constructed by chopping
off cusps of a finite volume complete locally symmetric manifold
of non-constant sectional curvature, and then doubling it along 
the boundary. Because the simplicial volume is additive under gluings
along boundary incompressible components with amenable
fundamental groups (see Kuessner~\cite{Kue-mul}) the double has nonzero
simplicial
volume, yet it admits no locally--$CAT(0)$ metric because
its fundamental group contains a nonabelian nilpotent subgroup
coming from the cusp (see Bridson and Haefliger~\cite[page 439]{BH}).

This idea can be pushed much farther by noting that by the item (III)
in the abstract the fundamental group
of any finite aspherical complex embeds into $\pi_1(D)$ for some 
$D$ in $\mathcal {AP}$, so it is very easy to find examples
when $\pi_1(D)$ is not $CAT(0)$.
For example, following the idea of Mess~\cite{Mess},
one finds $D$ in $\mathcal {AP}$ such that $\pi_1(D)$
contains the Baumslag--Solitar group $BS(1,2)$, 
which in turn contains a subgroup of dyadic rationals,
which cannot happen for $CAT(0)$ groups.
Here is a more sophisticated example based on an observation of 
Weinberger~\cite[Section 13]{Davis-ex}.

\begin{cor}\label{cor: word prob}
For any $n\ge 5$ there exists a closed aspherical PL $n$--manifold
with nonzero simplicial volume, and non-elementary
relatively hyperbolic fundamental group that has unsolvable 
word problem. 
\end{cor}
\begin{proof} 
By work of Collins and Miller~\cite{ColMil}
there exists a finite aspherical $2$--complex $X_0$ whose
fundamental group has unsolvable word problem. 
By a theorem of Stallings, $X_0$ is homotopy equivalent
to a finite aspherical complex $X$ that 
embeds into $\R^n$ for each $n\ge 4$ (see Stallings~\cite{Sta}, and
Drani\v{s}nikov and Repov\v{s}~\cite{DR}).
So the reflection group trick produces a closed aspherical 
PL $n$--manifold $B$ that retracts onto $X$, and by 
\fullref{thm: main}
there exists a closed aspherical PL manifold $D$ of dimension $n{+}1$
with nonzero simplicial volume, and non-elementary
relatively hyperbolic fundamental group that retracts
onto $X$. If a group has a solvable word problem, so do
all of its finitely generated subgroups, so
the word problem of $\pi_1(D)$ is unsolvable.
\end{proof}

\begin{rmk}
Note that the word problem of a finitely presented group 
$H$ is solvable if $H$ is
$CAT(0)$ (see Bridson and Haefliger~\cite[page 441]{BH}) 
or hyperbolic, or automatic, or  
asynchronously automatic, 
or residually finite,
or linear over any commutative ring 
(see Miller~\cite{Mil} and references therein).
%
Thus the group in 
\fullref{cor: word prob} satisfies none of these
properties.  
\end{rmk}

\section{Applications to assembly maps}

To motivate the discussion of this section we mention,
following Davis~\cite{JDav-nov}, 
some outstanding conjectures on the assembly maps
in $L$--theory and $K$--theory
\[
A^L_\pi\co H_\ast(B\pi;\mathbf L)\to L_\ast(\Z\pi)\qquad\quad
A^K_\pi\co H_\ast(B\pi;\mathbf K)\to K_\ast(\Z\pi)
\]
which are still open even when the Eilenberg--MacLane
space $B\pi$ can be realized as a closed aspherical manifold.
\smallskip

\begin{nhsc}
\label{nhsc}
For any group $\pi$,
the rational assembly maps $A^L_\pi\otimes\mathrm{id}_{\mathbb Q}$ and
$A^K_\pi\otimes\mathrm{id}_{\mathbb Q}$ are injective.
\end{nhsc}

\begin{amc}
\label{amc}
If $\pi$ is torsion-free,
the assembly maps $A^L_\pi$ and $A^K_\pi$ are isomorphisms.
\end{amc}

It is known that the assembly map in algebraic $K$--theory
is always injective in degrees $\ast <2$ with
cokernels $K_\ast(\Z\pi)$, $\tilde K_0(\Z\pi)$, $\mathrm{Wh}(\pi)$ 
for $\ast<0$, $\ast=0$, 
$\ast=1$, respectively. Thus the validity of the \fullref{amc}
in $K$--theory yields significant topological information. 
Other topological
consequences can be obtained via the surgery exact sequence,
in which the surgery obstruction map is closely related to the $L$--theory
assembly map.
One of the key conjectures in geometric topology is the following.\smallskip

\begin{rbc}
\label{rbc}
If $M$ is a compact aspherical
manifold with possibly nonempty boundary, then any homotopy
equivalence of manifold pairs $(N,\d N)\to (M,\d M)$ that 
restricts to a homeomorphism $\d N\to\d M$ must be homotopic, 
relative to the boundary, to a homeomorphism $N\to M$.
\end{rbc}

It is known that if $M$ is a compact
aspherical manifold of dimension $\ge 5$, then  
the \fullref{rbc} for $M$
holds if and only if $\mathrm{Wh}(\pi_1(M))$ vanishes
and the \fullref{amc} in $L$--theory holds for 
$\pi_1(M)$.

For the purposes of this paper, the key property of the assembly
map is naturality 
(which can be deduced eg from Hambleton and Pedersen~\cite{HamPed}).
More precisely, the $K$--theory assembly map
is natural with respect to any group homomorphism $\phi\co G\to\pi$,
while the $L$--theory assembly is only natural under the 
twisting preserving homomorphisms. 
Namely, the definition of $L(\Z\pi)$ involves a twisting, 
ie a homomorphism $w(\pi)\co\pi\to\Z_2$, and we require that
$w(G)=w(\pi)\circ\phi$.  
If $\pi$ is the fundamental group 
of a manifold $M$, then there is a natural twisting coming 
from the first Stiefel--Whitney class of $M$.

Suppose we have a retraction of groups $G\to \pi$, 
eg the retraction of the fundamental groups induced by the 
retraction of aspherical spaces in (I), (II). 
To see the idea of what follows, note that since
the Whitehead group depend covariantly on the group, we get 
an induced retraction $\mathrm{Wh}(G)\to \mathrm{Wh}(\pi)$, 
in particular, if  
$\mathrm{Wh}(G)$ vanishes, then so does $\mathrm{Wh}(\pi)$.

More generally, let $A\co\mathbf{Q}\to\mathbf {P}$
be a natural transformation of two covariant
functors $\mathbf{Q}$, $\mathbf {P}$ defined on and taking values
in the category of groups and group homomorphisms.
Applying $A$ to the inclusion $\pi\to G$ and the retraction $G\to\pi$,
we get the following commutative diagram.
\[
\xymatrix{
\mathbf{Q}(\pi)\ar[r]\ar[d]^{A_\pi}&
\mathbf{Q}(G)\ar[r]\ar[d]^{A_G}&
\mathbf{Q}(\pi)\ar[d]^{A_\pi}
\\
\mathbf{P}(\pi)
\ar[r]&
\mathbf{P}(G)
\ar[r]& 
\mathbf{P}(\pi)
}
\]
Since the rows are induced by the retraction,
the composition of the two horizontal maps 
in each row is the identity, so the first map is injective,
and the second map is onto. 
An easy diagram chase shows that if $A_G$
is injective, then so is $A_\pi$, and moreover,
the same statement holds if the word ``injective'' 
is replaced by ``surjective'', or
``isomorphism'', or ``split injective''.
Again, to apply the above to the $L$--theory assembly map 
we have to assume that the retraction $G\to \pi$ preserves 
the twisting, but this assumption is not very restrictive
because eg given any twisting $w\co\pi\to\Z_2$, there is 
a twisting of $G$ that makes $G\to\pi$ twisting-preserving, 
namely, the twisting obtained by composing $G\to\pi$ and $w$.

I suspect that this idea has been 
well-known for a long time, eg 
since any closed aspherical $n$--manifold $M^n$ is a retract 
of a closed aspherical $(n{+}1)$--manifold, namely $M^n\times S^1$, 
one sees that if the \fullref{amc} holds 
for all closed aspherical $(n{+}1)$--manifolds, then it also holds
all closed aspherical $n$--manifolds. Here the retraction 
$M^n\times S^1\to M$ preserves the twisting given by 
the first Stiefel--Whitney class $w_1$, because it 
pulls back the stable tangent bundle. 

Davis used this idea in~\cite[Sections 11]{Davis-ex}
to show that the \fullref{amc} holds for
all closed aspherical manifolds, it also holds
of all groups with $B\pi$ a finite complex. Indeed, by the
reflection group trick any finite aspherical simplicial
complex $K$ is a retract of a closed aspherical manifold $M$.
Furthermore, given any twisting $w\co\pi_1(K)\to\Z_2$
one can even choose $M$ with $w_1(M)=w\circ r_\ast$, where
$r\co M\to K$ denotes the retraction.
Indeed, $w$ can be realized as the first Stiefel--Whitney class of
a compact manifold $F$, that is a thickening of $K$ (eg
embed $K$ into a Euclidean space, take a regular neighborhood, 
and then let $F$ to be the line bundle over the regular neighborhood
with the first Stiefel--Whitney class $w$).
Now use $F$ as the fundamental chamber in the reflection group trick
to produce a closed aspherical manifold $M$ that retracts onto $F$, 
and hence onto $K$.
By a result of Davis~\cite[Proposition 1.4]{Davis-some-asph-mfld},
the retraction $M\to F$ pulls back the stable tangent bundle, so it
preserves the twisting.

Finally, specializing to the setting of \fullref{thm: main},
we get, among other things,
the results on the \fullref{amc}
and the Novikov Conjecture stated in the introduction.
By \fullref{prop: retract w1} if $m=1$, then  
the retraction $r_1$ 
preserves the natural twistings given by $w_1$.

\section{Questions}

\begin{quest}\rm Are there closed aspherical manifolds with hyperbolic 
fundamental group which are not homotopy equivalent 
to a triangulated manifold? If you drop the requirement that the
fundamental group is hyperbolic, such examples were constructed
by Davis and Januszkiewicz~\cite{DJ}. 
\end{quest}

\begin{quest}\rm
Suppose we start from a finite volume real hyperbolic 
manifold with toral cusp cross-section, and fill them by 
a ``sufficiently complicated'' Dehn surgery (see Anderson~\cite{And-dehn}) so that the result
is a closed aspherical manifold $M$ of nonpositive sectional curvature.
Under what conditions $M$ has nonzero simplicial volume? 
(By \fullref{ex: mosher-sageev}
this is true if each cusp has a cross-section of injectivity 
radius $>\pi$).
In particular, do the aspherical homology 
$4$--spheres constructed by Ratcliffe and Tschantz~\cite{RatTsc}
have nonzero simplicial volume? Note that there are many homology 
$3$--spheres carrying hyperbolic metrics, so that their simplicial volumes
are nonzero.
\end{quest}

\section{Acknowledgements}

I gratefully acknowledge support by the 
NSF grants \# DMS--0503864, DMS--0352576, and am thankful
to James F Davis, Erik K Pedersen, and Andrew Ranicki
for helping me to locate the correct reference to the naturality
of the assembly map.

\bibliographystyle{gtart}
\bibliography{link}

\end{document}